\documentclass{article}
\pdfoutput=1
\usepackage{amssymb,amsmath}
\usepackage{amsfonts}
\usepackage{latexsym}
\newtheorem{Lem}{Lemma}
\newtheorem{Thm}{Theorem}
\newtheorem{Pro}{Proposition}
\newtheorem{Rem}{Remark}
\newtheorem{Cor}{Corollary}
\newtheorem{Exa}{Example}
\newtheorem{Def}{Definition}
\def\C{{\mathbb C}}
\def\R{{\mathbb R}}

\def\H{{\mathbb H}}
\def\M{{\mathbb M}} 
\def\S{{\mathbb S}}

\title{Quaternionic Regularity via Analytic Functional Calculus}
\author{Florian-Horia Vasilescu\\
\small Department of Mathematics, University of Lille,\\
\small 59655 Villeneuve d'Ascq, France\\
\small e-mail: florian.vasilescu@univ-lille.fr}
\date{}

\begin{document}

\maketitle

\begin{abstract} Denoting by $\M$ the complexification of the quaternionic algebra $\H$, we characterize the family of those
$\M$-valued functions, defined on subsets of $\H$, whose values
are actually quaternions, using  an intrinsic approach. In particular, we show that the slice quaternionic regularity can be characterized
via a Cauchy type transform acting on the space of analytic 
$\M$-valued stem functions. 
\end{abstract}
\medskip



\section{Introduction}\label{I}

The quaternions form a unital non commutative division algebra, with numerous applications in mathematics and physics. 
In mathematics, one of the most important investigation in the quaternionic context has been to find a convenient manner to express the ''analyticity`` of functions depending on quaternions.  Among the pioneer contributions in this direction one should mention the works \cite{MoTh} and \cite{Fue}. Concerning the connections with physics, the foundations of the quaternion quantum mechanics can be found in the work \cite{FiJaScSp}. 

In 2007, a concept of {\it slice regularity} for functions of one quaternionic variable was introduced  in \cite{GeSt}, leading to a large development sythesized in \cite{CoSaSt} (see also 
\cite{GhPe}, \cite{GhMoPe}, \cite{CoGaKi} etc.).

Unlike in \cite{GeSt}, in the paper \cite{Vas6} the regularity
of a quaternionic-valued function was investigated via the analytic functional calculus acting on quaternions.  This was done by considering the algebra of quaternions as a real subalgebra of the complex algebra of  $2\times 2$ matrices with complex entries. This well known matrix representation allowed us  
to view  the quaternions as linear operators  on a complex space, and thus commuting with the complex numbers. Moreover, each quaternion was regarded as a normal operator, having a spectrum which was used to define various compatible functional calculi, including the analytic one. Specifically, 
 a ''quaternionic regular function`` was obtained by a pointwise construction of  the analytic functional calculus with stem functions on a conjugate symmetric open set $U$ in the complex plane, applied to quaternions whose spectra was in $U$,  via the  matrix version of  Cauchy's formula, with no need of slice derivatives.
In fact,  all of the  ''regular functions``, regarded as   quaternionic Cauchy transforms of stem functions, happened to have
some unexpected multiplicative properties, and several properties of slice regular functions (see \cite{CoSaSt})  was recaptured.

In the present paper we treat similar problems to those from the first part of \cite{Vas6}, in the  context of the abstract Hamilton's algebra, embedding the real $C^*$-algebra of quaternions into its complexification, organized as a complex $C^*$-algebra. Using elements of spectral theory, we succed to considerably simplify the proof the main corresponding results from \cite{Vas6}. Nevertheless, the second part of 
\cite{Vas6}, dealing with real and quaternionic operators, is not covered by the present text. Unlike in \cite{Vas6}, the actual arguments are not only much simpler but the framework is intrisic, that is, it does not depend of any representation of Hamilton's algebra.
 
Let us briefly describe the contents of this work. The next chapter is dedicated to some preliminaries, including Hamilton's algebra of quaternions, denoted by $\H$, and its complexification $\M$, endowed with a unique $C^*$-algebra structure. The class of $\M$-valued slice regular functions is also mentioned.  

In the third chapter, the spectrum of a quaternion as an
element of the algebra $\M$ is introduced (see Remark \ref{spectrum}), and  the slice regularity of an associated 
 $\M$-valued Cauchy kernel is proved. As the spectrum of a quaternion consists of at most two points, the structure of the associated spectral projections is also exhibited, for later use.  

The representation of a quaternion in terms of a resolution of 
identity, as a particular case of the concept of scalar operator
(see \cite{DuSc}, Part III), allows us the define a functional calculus with arbitrary $\M$-valued functions for each quaternion,
naturally extended to sets of quaternions in a pointwise manner
(see Definition \ref{funccalc}). In addition, the function given by
formula (\ref{GFC2}) is $\H$-valued if and only if it is constructed 
by using a stem function, via  Theorem \ref{sym_spec}, leading to 
a general functional calculus with arbitrary functions (see 
Theorem \ref{gfc}). 
 
A quaternionic Cauchy transform is introduced in the fifth chapter, which can be defined for all analytic $\M$-valued functions, but
especially of interest when working with analytic stem functions. 

In the last chapter the analytic functional calculus for quaternionic functions is obtained (see Theorem \ref{H_afc}), as a particular case of the general functional calculus given by Theorem \ref{gfc}. We also recapture one of  the main results from \cite{Vas6}, showing   that the slice regularity of a quaterninic function is equivalent with its property of being the Cauchy transform of a stem function (Theorem \ref{equiv-ons-dom1}). The main ingredients of this result are the 
representation formula (\ref{reprez}) and a similar result from
\cite{CoSaSt} (Lemma 4.3.8).

\section{Preliminaries} 
 
\subsection{Hamilton's Algebra}\label{HA}

For the sake of completeness, and to fix the notation and terminology, we start this discussion with some well known facts. Abstract Hamilton's algebra $\H$ is the four-dimensional $\R$-algebra with 
unit $1$, generated by the ''imaginary units`` $\{\bf{j,k,l}\}$,  
which satisfy
$$
{\bf jk=-kj=l,\,kl=-lk=j,\,lj=-jl=k,\,jj=kk=ll}=-1.
$$

We may assume that $\H\supset\R$ identifying every  number
$x\in\R$ with the element $x1\in\H$.

The algebra $\H$ has a natural multiplicative norm given by
$$
\Vert {\bf x}\Vert=\sqrt{x_0^2+x_1^2+x_2^2+x_0^2},\,\,{\bf x}= x_0+x_1{\bf j}+x_2{\bf k}+x_3{\bf l},\,\,x_0,x_1,x_2,x_3\in\R,
$$
and a natural involution
$$
\H\ni{\bf x} = x_0+x_1{\bf j}+x_2{\bf k}+x_3{\bf l}\mapsto
{\bf x}^*= x_0-x_1{\bf j}-x_2{\bf k}-x_3{\bf l}\in\H.
$$

Every element ${\bf x}\in\H\setminus\{0\}$ is invertible, and ${\bf x}^{-1}=
\Vert {\bf x}\Vert^{-2}{\bf x}^*$.

For an arbitrary quaternion ${\bf x}= x_0+x_1{\bf j}+x_2{\bf k}+x_3{\bf l},\,\,x_0,x_1,x_2,x_3\in\R$, we set $\Re{\bf x}=x_0=
({\bf x}+{\bf x}^*)/2$, and $\Im{\bf x}=x_1{\bf j}+x_2{\bf k}+x_3{\bf l}=({\bf x}-{\bf x}^*)/2$, that is, the {\it real} and 
the {\it imaginary part} of ${\bf x}$, respectively.

Using one possible equivalent definition, a real $C^*$-algebra is a real Banach $*$-algebra
$A$ satisfying the $C^*$-identity $\Vert a^*a\Vert=\Vert a\Vert^2$
for all $a\in A$, also having the property  $\Vert a^*a\Vert\le  \Vert a^*a+b^*b\Vert$ for all $a,b\in A$ (see for instance \cite{Ros}). It is clear that the algebra $\H$ is a real  $C^*$-algebra.

\subsection{The Associated Complex $C^*$-Algebra}

The ''imaginary units`` $\bf j,k,l$ of the algebra $\H$ will be considered independent of the  imaginary unit $i$ of the complex plane $\C$. Specifically, we construct the complexification $\C\otimes_\R\H$
of the $\R$-vector space $\H$ (see also \cite{GhMoPe}), which will be
identified with the  direct sum $\M=\H+i\H$, with the  the natural multiplicative structure given by
$$
({\bf x}_1 + i{\bf x}_2)({\bf y}_1 + i{\bf y}_2)=({\bf x}_1{\bf y}_1-{\bf x}_2{\bf x}_2) + i({\bf x}_1{\bf y}_2+
{\bf x}_2{\bf y}_1),\,\,{\bf x}_1,{\bf x}_2,{\bf y}_1,{\bf y}_2\in\H, 
$$
Of course, the algebra $\M$ contains the complex field $\C$, that is,  
every complex number $z=s+it,\,s,t\in\R$ is identified with the element $s1+it1\in\M$. In this way, $\M$ becomes an associative complex algebra, with unit $1$ and 
involution $({\bf x}_1+i{\bf x}_2)^*={\bf x}_1^*-i{\bf x}_2^*$, where ${\bf x}_1,{\bf x}_2\in\H$ are arbitrary, making $\M$ an involutive
algebra. Moreover,
in the algebra $\M$, the elements of $\H$ commute with all complex numbers.  

In the algebra $\M$ there also exists a natural conjugation given
by $\bar{\bf a}={\bf b}-i{\bf c}$, where ${\bf a}={\bf b}+i{\bf c}$ is arbitrary in $\M$, with ${\bf b},{\bf c}\in\H$ (see also
\cite{GhMoPe}). Note that $\overline{\bf a+b}=\bar{\bf a}+\bar{\bf b}$, and  $\overline{\bf ab}=\bar{\bf a}\bar{\bf b}$, in particular  $\overline{r\bf a}=r\bar{\bf a}$ for all ${\bf a},{\bf b}\in\M$, and $r\in\R$.  Moreover, $\bar{{\bf a}}={\bf a}$ if and only if ${\bf a}\in\H$, which is a useful characterization of the elements of $\H$ among 
those of $\M$. 

Using some results from \cite{Pal} (see also \cite{Ros}), the  algebra $\M$ may be endowed with a unique $C^*$-algebra structure, containing the algebra $\H$ as a real $C^*$-algebra, via the natural enbedding. 

In our particular case, we can apply a more direct procedure, via a standard matricial representation (as in \cite{Vas6}). Namely, we have the following.

\begin{Thm}\label{HinM} The complex algebra $\M$ has a unique 
$C^*$-algebra structure such that $\H$ is a real $C^*$-subalgebra
of $\M$. 
\end{Thm}

{\it Sketch of proof.} There exists a $*$-isomorphism between the 
involutive algebra $\M$ and the algebra $\M_2$ consisting of all $2\times2$-matrices with complex entries, which is a complex $C^*$-algebra. Then the norm of $\M_2$ induces a norm on $\M$, making it a $C^*$-subalgebra. We omit the details.   

\subsection{Slice Regular Functions}

There exists a large literature dedicated to a concept of "slice regularity``, which is a form of holomorphy in the context of 
quaternions (see for instance \cite{CoSaSt} and the works quoted 
within).  

For $\M$-valued functions defined on subsets of $\H$, the concept of slice regularity (see \cite{CoSaSt}) is defined as follows. 

Let $\S=\{\mathfrak{s}=
x_1{\bf j}+
x_2{\bf k}+x_3{\bf l};x_1,x_2,x_3\in\R,x_1^2+x_2^2+x_3^2=1\}$, that is, the unit sphere of purely imaginary quaternions. 
It is clear that $\mathfrak{s}^*=-\mathfrak{s}$, and so 
$\mathfrak{s}^2=-1,\,\mathfrak{s}^{-1}=-\mathfrak{s}$, and
$\Vert \mathfrak{s}\Vert=1$ for all $\mathfrak{s}\in\S$.

Let also $\Omega\in\H$ be an open set, and let $F:\Omega\mapsto\M$ be a differentiable function. In the spirit of \cite{CoSaSt}, we say that $F$ is {\it right slice regular} on $\Omega$ if for all 
$\mathfrak{s}\in\S$,
$$
 \bar{\partial}_\mathfrak{s}F(x+y\mathfrak{s}):
 =\frac{1}{2}
\left(\frac{\partial}{\partial x}+R_\mathfrak{s}\frac{\partial}{\partial y}\right)F(x+y\mathfrak{s})=0,
$$
on the set $\Omega\cap(\R+\R\mathfrak{s})$,
where $R_\mathfrak{s}$ is the right multiplication of the elements of $\M$ by $\mathfrak{s}$. 

Note that, unlike in \cite{GeSt}, we use the right slice regularity rather than the left one because of a reason to be later seen.
Nevertheless, a left slice regularity can also be defined via the left multiplication of the elements of $\M$ by elements from 
$\S$. In what follows, the right slice regularity will be simply called {\it slice regularity}. 

We are particularly interested by slice regularity of 
$\H$-valued functions, but the concept is valid for $\M$-valued
functions and plays an important role in our discussion.

\begin{Exa}\label{exa_sr}\rm 
(1) The convergent series of the form $\sum_{k\ge 0}a_k{\bf q}^k$, 
on quaternionic balls 
$\{{\bf q}\in\H;\Vert {\bf q}\Vert<r\}$, with $r>0$ and $a_k\in\H$ for all $k\ge0$, are $\H$-valued slice regular on their domain of definition. 
In fact, if actually $a_k\in\M$, such functions are $\M$-valued  slice regular on their domain of definition.

(2) An important example of slice regular $\M$-valued function
will be further given by Example \ref{SliceregCk}.
\end{Exa}

\section{Spectrum of a Quaternion}

In the complex algebra $\M$ we have a natural concept of spectrum,
which can be easily described in the case of quaternions. In fact,
this spectrum coincides with that one introduced in \cite{Vas6}\,(see also \cite{CoSaSt}).


\begin{Rem}\label{spectrum}\rm (1) As each quaternion commutes in $\M$ with every complex number,
we have the identities
$$
(\lambda-{\bf x}^*)(\lambda-{\bf x})=(\lambda-{\bf x})(\lambda-{\bf x}^*)=\lambda^2-
\lambda({\bf x}+{\bf x}^*)+\Vert {\bf x}\Vert^2\in\C,
$$
for all $\lambda\in\C$ and ${\bf x}\in\H$. Therefore, the element
$\lambda-{\bf x}\in\M$ is invertible if and only if the complex number $\lambda^2-2\lambda\Re{\bf x}+\Vert {\bf x}\Vert^2$ is 
nonnull, and in that case 
$$
(\lambda-{\bf x})^{-1}=\frac{1}{\lambda^2-2\lambda\Re{\bf x}+\Vert {\bf x}\Vert^2}(\lambda-{\bf x}^*).
$$ 
Hence, the element      
$\lambda-{\bf x}\in\M$ is not  invertible if and only if $\lambda=
\Re{\bf x}\pm i\Vert\Im{\bf x}\Vert$. In this way, the {\it spectrum} of a quaternion ${\bf x}\in\H$ is given by the equality  $\sigma({\bf x})=\{s_\pm(\bf x)\}$, where 
$s_\pm(\bf x)=\Re{\bf x}\pm i\Vert\Im{\bf x}\Vert$ are the {\it eigenvalues} of $\bf x$.   

(2) As usually, the {\it resolvent set} $\rho({\bf x})$ of a quaternion ${\bf x}\in\H$ is the set $\C\setminus\sigma({\bf x})$, while the function
 $$\rho({\bf x})\ni\lambda\mapsto(\lambda-{\bf x})^{-1}\in\M$$
 is the {\it resolvent $($function$)$} of ${\bf x}$, which is a 
$\M$-valued analytic function on $\rho({\bf x})$. 

(3)  Note that two quaternions ${\bf x},
{\bf y}\in\H$ have the same spectrum if and only if  $\Re{\bf x}=
\Re{\bf y}$ and $\Vert\Im{\bf x}\Vert=\Vert\Im{\bf y}\Vert$.

(4) As before, let $\S$ be the unit sphere of purely imaginary quaternions. It is  clear that  every quaternion 
${\bf q}\in \H\setminus\R$ can be  written as ${\bf q} = x + y\mathfrak{s}$, where $x,y$ are real numbers, with $x=\Re{\bf q}$,
$y=\pm \Vert \Im{\bf q}\Vert$, and $\mathfrak{s}=\pm\Im{\bf q}/\Vert \Im{\bf q}\Vert\in\S$. Anyway, we always have 
$\sigma({\bf q})=\{x\pm iy\}$, because $\Im{\bf q}=y\mathfrak{s}$.
Note that, for fixed real numbers $x,y$, the spectrum of $\bf q$ does not depend on $\mathfrak{s}$. Thus, for every $\lambda=u+iv\in\C$ with $u,v\in\R$, we have $\sigma(u+v\mathfrak{s})=\{\lambda,\bar{\lambda}\}$. 
 
 (5) Fixing an element $\mathfrak{s}\in\S$, 
we define an isometric $\R$-linear map from the complex plane $\C$ into the  algebra $\H$, say $\tau_\mathfrak{s}$,  given by $\tau_\mathfrak{s}(u+iv)=u+v\mathfrak{s},\,u,v\in\R$. For every subset $A\in\C$, we put
\begin{equation}\label{embedA}
A_\mathfrak{s}=\{x+y\mathfrak{s}; x,y\in\R,x+iy\in A\}=\tau_\mathfrak{s}(A).
\end{equation}
Note that, if $A$ is open in $\C$, then $A_\mathfrak{s}$ is open 
in the $\R$-vector space $\C_\mathfrak{s}$.
 \end{Rem}

\begin{Def}\label{MCauchy}\rm  The $\M$-{\it valued Cauchy kernel} on the open set $\Omega\subset\H$ is given by
\begin{equation}\label{MvalCk}
\rho({\bf q})\times\Omega\ni(\zeta,{\bf q})\mapsto (\zeta-{\bf q})^{-1}\in\M.
\end{equation}
\end{Def}

\begin{Exa}\label{SliceregCk}\rm The $\M$-{\it valued Cauchy kernel} on the open set $\Omega\subset\H$ is slice regular. 
Specifically, choosing an arbitrary relatively open set $V\subset\Omega\cap(\R+\R\mathfrak{s})$, and fixing 
$\zeta\in\cap_{{\bf q}\in V}\rho({\bf q})$, we can write for 
${\bf q}\in V$ the equalities
$$
\frac{\partial}{\partial x}(\zeta-x-y\mathfrak{s})^{-1}=
(\zeta-x-y\mathfrak{s})^{-2},
$$
$$
R_\mathfrak{s}\frac{\partial}{\partial y}(\zeta-x-
y\mathfrak{s})^{-1}=-(\zeta-x-y\mathfrak{s})^{-2},
$$
because $\zeta$,  $\mathfrak{s}$ and 
$(\zeta-x-y\mathfrak{s})^{-1}$ commute in $\M$. Therefore,
$$
\bar{\partial}_\mathfrak{s}((\zeta-q)^{-1})=\bar{\partial}_\mathfrak{s}((\zeta-x-y\mathfrak{s})^{-1})=0,
$$ 
implying the assertion. 
\end{Exa}


\begin{Rem}\label{eigenval}\rm (1) The discussion about the spectrum
of a quaternion can be enlarged, keeping the same background. 
Specifically, we may regard an element $\bf q\in\H$ as a left multiplication operator on the $C^*$-algebra $\M$, denoted by $L_{\bf q}$, and given by $L_{\bf q}{\bf a}={\bf qa}$ for all 
${\bf a}\in\M$. It is easily seen that $\sigma(L_{\bf q})=\sigma(\bf q)$. In this context, we may find the eigenvectors of 
$L_{\bf q}$, which would be of interest in what follows.  Therefore, we should look for solutions of the equation
${\bf q}\nu=s\nu$ in the algebra $\M$, with $s\in\sigma(\bf q)$. 
Writing ${\bf q}=q_0+\Im{\bf q}$  with $q_0\in\R$, $s_\pm=q_0\pm i\Vert\Im{\bf q}\Vert$, and  
$\nu={\bf x}+i{\bf y}$ with ${\bf x},{\bf y}\in\H$, we obtain 
the equivalent equations
$$
(\Im{\bf q}){\bf x}=\mp\Vert\Im{\bf q})\Vert{\bf y},\,\, 
(\Im{\bf q}){\bf y}=\pm\Vert\Im{\bf q})\Vert{\bf x},
$$
leading to the solutions
$$
\nu_\pm({\bf q})=\left(1\mp i\frac{\Im\bf q}{\Vert{\Im \bf q}\Vert}\right){\bf x}
$$
of the equation ${\bf q}\nu_\pm=s_\pm\nu_\pm$, where 
${\bf x}\in\H$ is arbitrary, provided 
${\Im \bf q}\neq 0$.  

When ${\Im \bf q}=0$, the  solutions are given by  
$\nu=q_0{\bf a}$, with ${\bf a}\in\M$ arbitrary. 

(2) Every quaternion $\mathfrak{s}\in\mathbb{S}$ may be associated 
with two elements $\iota_\pm(\mathfrak{s})=(1\mp i\mathfrak{s})/2$ in $\M$, which are commuting idempotents
such that $\iota_+(\mathfrak{s})+\iota_-(\mathfrak{s})=1$ and
$\iota_+(\mathfrak{s})\iota_-(\mathfrak{s})=0$. For this reason,
setting $\M^\mathfrak{s}_\pm=\iota_\pm(\mathfrak{s})\H$, we have 
a direct sum decomposition $\M=\M^{\mathfrak{s}}_+ +
\M^\mathfrak{s}_-$. Indeed, it is clear that $\M^\mathfrak{s}_+\cap\M^\mathfrak{s}_-=\{0\}$. In addition, if ${\bf a}={\bf u}+i{\bf v}$,
with $\bf u,v\in\H$, the equation $\iota_+(\mathfrak{s}){\bf x}+
\iota_-(\mathfrak{s}){\bf y}=\bf a$
has the solution ${\bf x}={\bf u+\mathfrak{s}v},\,
{\bf y}={\bf u-\mathfrak{s}v}\in\H$, because  
$\mathfrak{s}^{-1}=-\mathfrak{s}$. Therefore, we also have, $\M=\M^\mathfrak{s}_+ +\M^\mathfrak{s}_-$.

In particular, if ${\bf q}\in\H$ and  ${\Im \bf q}\neq 0$, setting 
$\mathfrak{s}_{\tilde{\bf q}}=\tilde{\bf q}\Vert\tilde{\bf q}\Vert
^{-1 }$, where $\tilde{\bf q}=\Im\bf q$,
the elements $\iota_\pm(\mathfrak{s}_{\tilde{\bf q}})$ are idempotents, as above. Moreover,
\begin{equation}\label{eigenvalue_f}
{\bf q}\iota_\pm(\mathfrak{s}_{\tilde{\bf q}}){\bf h}=s_\pm({\bf q})\iota_\pm(\mathfrak{s}_{\tilde{\bf q}}){\bf h},\,\,{\bf h}\in\H.
\end{equation}
\end{Rem}

The next result provides explicit formulas of the spectral 
projections (see \cite{DuSc}, Part I, Section VII.1) associated to the operator $L_{\bf q},\,{\bf q}\in\H.$ 
Of course, this is not trivial only if ${\bf q}\in\H\setminus\R$
because if ${\bf q}\in\R$, its spectrum a real singleton, and the 
the only spectral projection is the identity.


\begin{Lem}\label{spproj} Let ${\bf q}\in\H\setminus\R$ be fixed. The spectral projections associated to $s_\pm({\bf q})$ are given by 
$$
P_\pm({\bf q}){\bf a}=\iota_\pm(\mathfrak{s}_{\tilde{\bf q}}){\bf a},\,\,{\bf a}\in\M.
$$
Moreover, $P_+({\bf q})P_-({\bf q})=P_-({\bf q})P_+({\bf q})=0$, and 
$P_+({\bf q})+P_-({\bf q})$ is the identity on $\M$.

When ${\bf q}\in\R$, the corresponding spectral projection is the identity on $\M$. 
\end{Lem}

{\it Proof.} Let us fix a quaternion $\bf q$ with ${\Im\bf q}\neq 0$. Next,  write the general formulas for its spectral projections. Setting $s_\pm= 
s_\pm({\bf q})$, the points $s_+,s_-$ are distinct and  not real. We fix an $r>0$ sufficiently small such that,
setting $D_\pm:=\{\zeta\in \rho({\bf q});\vert\zeta-s_\pm\vert\le r\}$, we have 
$D_\pm\subset \rho({\bf q})$ and $D_+\cap D_-=\emptyset$. Then 
we have
$$
P_\pm({\bf q})=\frac{1}{2\pi i}\int_{\Gamma_\pm} (\zeta-L_{\bf q})^{-1}d\zeta
$$
where $\Gamma_\pm$ is the boundary of $D_\pm$.  

Using the equality $L_{\bf q}\nu_\pm({\bf q})=
s_\pm({\bf q})\nu_\pm({\bf q})$
(see Remark \ref{eigenval}), for every $\zeta\in\rho({\bf q})$ and ${\bf h}\in\H$, we have 
$$
(\zeta-L_{\bf q})^{-1}(1\mp i\mathfrak{s}_{\tilde{\bf q}}){\bf h}=(\zeta- s_\pm)^{-1}(1\mp  i\mathfrak{s}_{\tilde{\bf q}}){\bf h},
$$
by Remark \ref{eigenval}. Therefore,
$$
P_+({\bf q})(1\mp i\mathfrak{s}_{\tilde{\bf q}}){\bf h}=\frac{1}{2\pi i}\int_{\Gamma_+} (\zeta- s_\pm)^{-1}(1\mp i\mathfrak{s}_{\tilde{\bf q}}){\bf h} d\zeta, 
$$
and
$$
P_-({\bf q})(1\mp i\mathfrak{s}_{\tilde{\bf q}}){\bf h}=\frac{1}{2\pi i}\int_{\Gamma_-} (\zeta- s_\pm)^{-1}(1\mp i\mathfrak{s}_{\tilde{\bf q}}){\bf h} d\zeta.
$$
Using Cauchy's formula, we deduce that
$$
P_+({\bf q})(1-i\mathfrak{s}_{\tilde{\bf q}}){\bf h}=(1-i\mathfrak{s}_{\tilde{\bf q}}){\bf h},\,P_+({\bf q})(1+i\mathfrak{s}_{\tilde{\bf q}}){\bf h}=0,
 $$
and
$$
P_-({\bf q})(1-i\mathfrak{s}_{\tilde{\bf q}}){\bf h}=0,\,P_-({\bf q})(1+i\mathfrak{s}_{\tilde{\bf q}}){\bf h}=(1+i\mathfrak{s}_{\tilde{\bf q}}){\bf h},
 $$
for all $\bf h\in\H$. 

Fixing an arbitrary element ${\bf a}={\bf u}+i{\bf v}\in\M$, writing 
$$
{\bf a}=\iota_+(\mathfrak{s}_{\tilde{\bf q}})({\bf u}+\mathfrak{s}_{\tilde{\bf q}}{\bf v})+\iota_-(\mathfrak{s}_{\tilde{\bf q}})({\bf u}-\mathfrak{s}_{\tilde{\bf q}}{\bf v})
$$
with ${\bf u}\pm\mathfrak{s}_{\tilde{\bf q}}{\bf v}\in\H$ (see Remark  \ref{eigenval}(2)), and  noticing that $\iota_\pm(\mathfrak{s}_{\tilde{\bf q}})\mathfrak{s}_{\tilde{\bf q}}=\pm i\iota_\pm(\mathfrak{s}_{\tilde{\bf q}})$, as $P_\pm({\bf q})$ are $\C$-linear, we  obtain

\begin{equation}\label{formproj}
P_\pm({\bf q}){\bf a}=\iota_\pm(\mathfrak{s}_{\tilde{\bf q}}){\bf a},\,\,{\bf a}\in\M,
\end{equation} 
 which are precisely the formulas from the statement.

The properties $P_+({\bf q})P_-({\bf q})=P_-({\bf q})P_+({\bf q})=0$, and 
$P_+({\bf q})+P_-({\bf q})$ is the identity on $\M$ are direct 
consequences of the analytic functional calculus associated to a
fixed element ${\bf q}\in\H$ in the algebra $\M$. 
 
\medskip

By a slight abuse of terminology, the projections $P_\pm({\bf q})$
will be also called the {\it spectral projections} of $\bf q$. 
In fact, as formula (\ref{formproj}) shows, they depend only 
on the imaginary part of ${\bf q}$.

\section{A General Functional Calculus}
 
 In this section, starting from some spaces of $\M$-valued functions, defined on  subsets of the complex plane, we construct and characterize functional calculi, taking values in the quaternionic algebra $\H$. 
 
 
 \begin{Rem}\label{GFC}\rm Regarding the $C^*$-algebra $\M$ as a (complex) Banach space, and denoting by $\mathcal{B}(\M)$ the Banach space of all linear 
operators acting on $\M$, the operator $L_{\bf q},\,\bf q\in\H$, (see Remark \ref{eigenval}(1)) is a (very) particular case of a {\it scalar type} operator, as defined in \cite{DuSc}, Part III, XV.4.1. Its resolution of the identity consists of four projections $\{0,P_\pm({\bf q}),$I$\}$, including the null operator $0$ and the identity I on $\M$, where $P_\pm({\bf q})$ are the spectral projections of
$L_{\bf q}$, and its integral representation is  given by
$$
L_{\bf q}=s_+({\bf q})P_+({\bf q})+s_-({\bf q})P_-({\bf q})\in\mathcal{B}(\M), 
$$
via formulas (\ref{eigenvalue_f}) and (\ref{formproj}).
For every function $f:\sigma({\bf q})\mapsto\C$ we may define the operator
$$
f(L_{\bf q})=f(s_+({\bf q}))P_+({\bf q})+f(s_-({\bf q}))P_-({\bf q})\in\mathcal{B}(\M).
$$
which provides a functional calculus with arbitrary functions on 
the spectrum. More generally, we may extend this formula 
to functions of the form  $F:\sigma({\bf q})\mapsto\M$, putting
\begin{equation}\label{GFC1}
F(L_{\bf q})=F(s_+({\bf q}))P_+({\bf q})+F(s_-({\bf q}))P_-({\bf q}),
\end{equation}
and keeping this order, which is a ''left functional calculus``, not multiplicative, in general. It is this idea which leads  us to try to define  $\H$-valued functions on subsets of $\H$ via some $\M$-valued functions, defined on  subsets of $\C$.
\end{Rem}


\begin{Def}\label{consym}\rm (1) A subset $S\subset\C$ is said to be {\it conjugate symmetric} if $\zeta\in S$ if and only if 
$\bar{\zeta}\in S$. 

(2) A subset $A\subset\H$ is said to be {\it spectrally saturated} 
(see \cite{Vas6}) if whenever $\sigma({\bf h})=\sigma({\bf q})$ for some ${\bf h}\in\H$ and ${\bf q}\in A$, we also have ${\bf h}\in A$. 

For an arbitrary $A\subset\H$, we put $\mathfrak{S}(A)=
\cup_{{\bf q}\in A}\sigma({\bf q})$. 
We also  put $S_\H=\{{\bf q}\in\H;\sigma({\bf q})\subset S\}$ for an arbitrary subset $S\subset\C$.
\end{Def}

 
\begin{Rem}\label{consymr}\rm 
(1) If $A\subset\H$ is spectrally saturated, then $S=\mathfrak{S}(A)$ is conjugate symmetric, and conversely, if $S\subset\C$ is conjugate symmetric, then $S_\H$ is spectrally saturated, which
can be easily seen. Moreover, the assignment $S\mapsto S_\H$ is injective. Indeed, if $\lambda=u+iv\in S,\,u,v\in\R$, then $\lambda\in\sigma(u+v\mathfrak{s})$ for a fixed $\mathfrak{s}\in\S$. If 
$S_\H=T_\H$ for some $T\subset\C$, we must have ${\bf q}=u+v\mathfrak{s}\in T_\H$. Therefore $\sigma(u+v\mathfrak{s})
\subset T$, implying $\lambda\in T$, and so $S\subset T$. Clearly,
we also have $T\subset S$. 

Similarly, the assignment $A\mapsto\mathfrak{S}(A)$ is injective and  $A=S_\H$ if and only if $S=\mathfrak{S}(A)$. These two assertions are left to the reader.

(2)   If $\Omega\subset\H$ is an open spectraly saturated set, then $\mathfrak{S}(\Omega)\subset\C$ is open. To see that, let $\lambda_0=u_0+iv_0\in\mathfrak{S}(\Omega)$ be fixed, with $u_0,v_0\in\R$,  and let ${\bf q}_0=u_0+v_0\mathfrak{s}$, where 
$\mathfrak{s}\in\S$ is also fixed. Because $\Omega$ is spectrally 
saturated, we must have ${\bf q}_0\in\Omega$. Because the set $\Omega\cap\C_\mathfrak{s}$ is relatively open, there is a positive number
$r$ such that the open set 
$$
\{{\bf q}=u+v\mathfrak{s}; u,v\in\R,\vert{\bf q}-{\bf q}_0\vert<r\}
$$
is in $\Omega\cap\C_\mathfrak{s}$, where ${\bf q}=u+v\mathfrak{s}$. Therefore, the
set of the points $\lambda=u+iv$, satisfying $\vert\lambda-\lambda_0\vert<r$ is in $\mathfrak{S}(\Omega)$, implying that it is open. 

Conversely,
if $U\subset\C$ is open and conjugate symmetric,
the set $U_\H$ is also open via the upper semi-continuity of the spectrum (see \cite{DuSc}, Part I, Lemma VII.6.3.).

An important particular case is when $U={\mathbb D}_r:=\{\zeta\in\C;\vert\zeta\vert<r\}$, for some $r>0$. Indeed, if $\vert\kappa\vert
< r$ and  $\theta$ has the property $\sigma({\bf q})=\sigma({\bf h})$, from the equality $\{\Re(\kappa)\pm i\vert\Im(\kappa)\vert\}=
\{\Re({\bf h})\pm i\vert\Im({\bf h})\vert\}$ it follows that 
$\vert{\bf h}\vert<r$.

 (3) A subset $\Omega\subset\H$ is said to 
be {\it axially symmetric} if for every ${\bf q}_0=u_0+v_0\mathfrak{s}_0\in\Omega$ with $u_0,v_0\in\R$ and $\mathfrak{s}_0\in\S$, we also
have ${\bf q}=u_0+v_0\mathfrak{s}\in\Omega$ for all $\mathfrak{s}\in\S$. This concept is similar to the corresponding one in \cite{CoSaSt}, Definition 2.2.17. In fact, we have the following.


\begin{Lem}\label{UH_tildeU} A subset $\Omega\subset\H$ is axially
 symmetric if and only if it is spectrally saturated. 
\end{Lem}

{\it Proof.}  Let $\Omega\subset\H$ be axially
symmetric and ${\bf q}=u+v\mathfrak{s}\in\Omega$. Therefore, 
$\sigma({\bf q})=\{u\pm iv\}$. Let also ${\bf h}\in\H$ be such that
$\sigma({\bf h})=\sigma({\bf q})$. Hence ${\bf h}=u+w\mathfrak{s}'$.
Changing if necessaery, the sign of $\mathfrak{s}'$, we may assume, 
with no loss of generality, that $w=v$. Then ${\bf h}\in\Omega$,
showing that $\Omega$ is spectrally saturated.

Conversely, assuming that $\Omega$ is spectrally saturated, and fixing
an element ${\bf q}_0=u_0+v_0\mathfrak{s}_0\in\Omega$, then each element of the form ${\bf q}=u_0+v_0\mathfrak{s}$  has the same spectrum as ${\bf q}_0$, and thus it must belong to $\Omega$. 
Consequently, $\Omega$ is axially symmetric.

Nevertheless, we continue to use the expression  ''spectrally saturated set`` to designate an ''axially symmetric set``, because the
former is more compatible with our spectral approach.  
\end{Rem}

As noticed above, the algebra $\M$ is endowed with a conjugation given
by $\bar{\bf a}={\bf b}-i{\bf c}$, when ${\bf a}={\bf b}+i{\bf c}$, with ${\bf b},{\bf c}\in\H$. Note also that, because $\C$
is a subalgebra of $\M$, the conjugation of $\M$ restricted to
$\C$ is precisely the usual complex conjugation.
\medskip

The next definition has an old origin, going back to \cite{Fue}
(see for instance \cite{GhPe}).


\begin{Def}\label{stem}\rm Let $U\subset\C$ be conjugate symmetric, and let $F:U\mapsto\M$. We say that $F$ is a {\it stem function} if  $F(\bar{\lambda})=\overline{F(\lambda)}$ for all 
$\lambda\in U$. 
\end{Def}
\medskip

For an arbitrary conjugate symmetric subset $U\subset\C$, we put 
\begin{equation}
\mathcal{S}(U,\M)=\{F:U\mapsto\M;F(\bar{\zeta})=\overline{F(\zeta)},\zeta\in U\},
\end{equation}
that is, the $\R$-vector space of all $\M$-valued stem functions on $U$.  
Replacing $\M$ by $\C$, we denote by $\mathcal{S}(U)$ 
the real algebra of all $\C$-valued stem functions, which is an 
$\R$-subalgebra in $\mathcal{S}(U,\M)$. In addition, the space
$\mathcal{S}(U,\M)$ is a $\mathcal{S}(U)$-bimodule.


\begin{Def}\label{funccalc} Let  $U\subset\C$ be  conjugate symmetric.
For every $F:U\mapsto\M$ and all ${\bf q}\in U_\H$ we define 
a function $F_\H:U_\H\mapsto\M$, via the assignment 
\begin{equation}\label{GFC2}
U_\H\ni{\bf q}\mapsto F_\H({\bf q})=F(s_+({\bf q}))\iota_+(\mathfrak{s}_{\tilde{\bf q}})+F(s_-({\bf q}))\iota_-(\mathfrak{s}_{\tilde{\bf q}}) \in\M,
\end{equation} 
where $\tilde{\bf q}=\Im\bf q,\,\mathfrak{s}_{\tilde{\bf q}}=\tilde{\bf q}\Vert\tilde{\bf q}\Vert^{-1 }$, and  
$\iota_\pm(\mathfrak{s}_{\tilde{\bf q}})=2^{-1}
(1\mp i\mathfrak{s}_{\tilde{\bf q}})$.
\end{Def} 

Formula (\ref{GFC2}) is strongly related to formula  
(\ref{GFC1}) because the spectral projections $P_\pm({\bf q})$ are the left multiplications defined by $2^{-1}(1\mp i\mathfrak{s}_{\tilde{\bf q}})$ respectively, via formula (\ref{formproj}).

The next result is an intrinsic version of Theorem 1 from \cite{Vas6}, with a much shorter proof. 


\begin{Thm}\label{sym_spec} Let $U\subset\C$ be a conjugate symmetric subset, and let $F:U\mapsto\M$.
The element $F_\H({\bf q})$ is a quaternion for all ${\bf q}\in U_\H$ if and only if  $F$ is a stem function.
\end{Thm}

{\it Proof.}\,  We first assume that $F_\H({\bf q})$ is a quaternion for all ${\bf q}\in U_\H$. We fix a point $\zeta\in U$, supposing that $\Im\zeta>0$. Then we choose a quaternion 
${\bf q}\in U_\H$ with $\sigma({\bf q})=\{\zeta,\bar{\zeta}\}$. 
Therefore, $s_+({\bf q})=\zeta$ and $s_-({\bf q})=\bar{\zeta}$.
We write $F(\zeta)=F_{+1}+iF_{+2}, F(\bar{\zeta})=F_{-1}+iF_{-2}$, 
with $F_{+1},F_{+2}, F_{-1},F_{-2}\in\H$. 
According to (\ref{GFC2}), we infer that
$$
2F({\bf q})=F_{+1}+F_{+2}\mathfrak{s}_{\tilde{\bf q}}+F_{-1}-F_{-2}\mathfrak{s}_{\tilde{\bf q}}+i(-F_{+1}\mathfrak{s}_{\tilde{\bf q}}
+F_{+2}+F_{-2}+F_{-1}\mathfrak{s}_{\tilde{\bf q}}),
$$
so
$$
2\overline{F({\bf q})}=F_{+1}+F_{+2}\mathfrak{s}_{\tilde{\bf q}}+F_{-1}-F_{-2}\mathfrak{s}_{\tilde{\bf q}}+i(F_{+1}\mathfrak{s}_{\tilde{\bf q}}-F_{+2}-F_{-2}-F_{-1}\mathfrak{s}_{\tilde{\bf q}}).
$$
Because $F({\bf q})=\overline{F({\bf q})}$, we must have
$$
-F_{+1}\mathfrak{s}_{\tilde{\bf q}}
+F_{+2}+F_{-2}+F_{-1}\mathfrak{s}_{\tilde{\bf q}}=F_{+1}\mathfrak{s}_{\tilde{\bf q}}-F_{+2}-F_{-2}-F_{-1}\mathfrak{s}_{\tilde{\bf q}},
$$
which is equivalent to    
$$
F_{+2}+F_{-2}=(F_{+1}-F_{-1})\mathfrak{s}_{\tilde{\bf q}}
$$  
Assuming $F_{+1}\neq F_{-1}$, we deduce that 
$$
(F_{+1}-F_{-1})^{-1}(F_{+2}+F_{-2})=\mathfrak{s}_{\tilde{\bf q}}.
$$   
This equality is impossible because the left hand side depends
only on $\zeta$ and $\bar{\zeta}$ while the right hand side 
has infinitely many distinct values, when replacing $\bf q$ by another
element from the set  $\{{\bf h}\in\H;s_+({\bf h})=\zeta\}\subset U_\H$.  Therefore, the equality   
$F({\bf q})=\overline{F({\bf q})}$ implies the equalities
$F_{+1}=F_{-1}$ and $F_{+2}=-F_{-2}$, meaning that 
$\overline{F(\zeta)}=F(\bar{\zeta})$. 

If $\Im\zeta=0$, so $\zeta=x_0\in\R$, taking ${\bf q}=x_0$, we have $\sigma({\bf q})=\{x_0\}$, and $F({\bf q})=F(x_0)$ is a 
quaternion.

If $\Im\zeta<0$, applying the above argument to $\bar{\zeta}$ we 
obtain $\overline{F(\bar{\zeta})}=F({\zeta})$. Consequently, $F$
is a stem function on $U$. 

Conversely, if  $\overline{F(\zeta)}=F(\bar{\zeta})$ for all $\zeta\in U$, choosing a ${\bf q}\in\H$ with $\Im{\bf q}\neq0$,
and fixing $\zeta=s_+({\bf q})$, we obtain from (\ref{GFC2})
the equality
$$
2F_\H({\bf q})=F(\zeta)+F(\bar{\zeta})-i(F(\zeta)-F(\bar{\zeta})\mathfrak{s}_{\tilde{\bf q}}.
$$
Therefore,
$$
2\overline{F_\H({\bf q})}=F(\bar{\zeta})+F(\zeta)+i(F(\bar{\zeta})-F(\zeta))\mathfrak{s}_{\tilde{\bf q}},
$$
showing that  $F({\bf q})\in\H$ for all ${\bf q}\in U_\H$,
because the case ${\bf q}=x_0\in\R$ is evident.  

\begin{Cor} Let $U\subset\C$ be a conjugate symmetric subset, and let $f:U\mapsto\C$.
The element $f_\H({\bf q})$ is a quaternion for all ${\bf q}\in U_\H$ if and only if  $f$ is a stem function. 
\end{Cor}

\begin{Rem}\label{zeros}\rm   Let $U\subset\C$ be a conjugate symmetric set and let $F\in\mathcal{S}(U,\M)$ be arbitrary. We can easily describe the zeros of $F_\H$. 
Indeed, if $F_\H({\bf q})=F(s_+({\bf q}))\iota_+({\tilde{\bf q}})+F(s_-({\bf q}))\iota_-({\tilde{\bf q}})=0$, we must have 
$F(s_+({\bf q}))\iota_+({\tilde{\bf q}})=0$ and 
$F(s_-({\bf q}))\iota_-({\tilde{\bf q}})=0$, via a direct manipulation with the idempotents $\iota_\pm({\tilde{\bf q}})$. In other words, we must have $F(s_\pm({\bf q}))=\pm iF(s_\pm({\bf q}))\mathfrak{s}_{\tilde{\bf q}}$. Choosing another 
quaternion ${\bf h}$ with $s_+({\bf q})\in\sigma({\bf h})$ and
$\tilde{{\bf q}}\neq\tilde{{\bf h}}$, we obtain
$F(s_+({\bf q}))(\mathfrak{s}_{\tilde{\bf q}}-\mathfrak{s}_{\tilde{\bf h}})=0$. Therefore, $F(s_+({\bf q}))=0$ because 
$ \mathfrak{s}_{\tilde{\bf q}}-\mathfrak{s}_{\tilde{\bf h}}$ is
invertible. Similarly, $F(s_-({\bf q}))=0$. Conqequently, setting $\mathcal{Z}(F):=\{\lambda\in U;F(\lambda)=0\}$,
and $\mathcal{Z}(F_\H):=\{{\bf q}\in U_\H;F_\H({\bf q})=0\}$,
we must have
$$
\mathcal{Z}(F_\H)=\{{\bf q}\in U_\H;\sigma({\bf q})\subset\mathcal{Z}(F)\}. 
$$
\end{Rem}

For every subset  $\Omega\subset\H$, we denote by $\mathcal{F}(\Omega,\H)$ the set of all $\H$-valued functions on $\Omega$.

The next result offers an  {\it $\H$-valued general functional calculus} with arbitrary stem functions.

\begin{Thm}\label{gfc} Let $\Omega\subset\H$ be a spectrally
saturated set, and let $U=\mathfrak{S}(\Omega)$. The map
$$
{\mathcal S}(U,\M)\ni F\mapsto F_\H\in\mathcal{F}(\Omega,\H)
$$ 
is $\R$-linear, injective, and has the property $(Ff)_\H=F_\H f_\H$ for all
$F\in{\mathcal S}(U,\M)$ and $f\in{\mathcal S}(U)$. Moreover, the
restricted map 
$$
{\mathcal S}(U)\ni f\mapsto f_\H\in\mathcal{F}(\Omega,\H)
$$ 
is unital and multiplicative.
\end{Thm}

{\it Proof.}\, The map $F\mapsto F_\H$ is clearly $\R$-linear.
The injectivity of this map follows from Remark \ref{zeros}. 
Note also that
$$
F_\H({\bf q}) f_\H({\bf q})=(F(s_+({\bf q}))\iota_+(\mathfrak{s}_{\tilde{\bf q}})+F(s_-({\bf q}))\iota_-(\mathfrak{s}_{\tilde{\bf q}}))\times 
$$
$$
(f(s_+({\bf q}))\iota_+(\mathfrak{s}_{\tilde{\bf q}})+f(s_-({\bf q}))\iota_-(\mathfrak{s}_{\tilde{\bf q}})= 
$$
$$
(Ff)(s_+({\bf q}))\iota_+(\mathfrak{s}_{\tilde{\bf q}})+(Ff)(s_-({\bf q}))\iota_-(\mathfrak{s}_{\tilde{\bf q}})=(Ff)_\H ({\bf q}),
$$ 
because $f$ is complex valued, and by the properties of the idempotents $\iota_\pm(\mathfrak{s}_{\tilde{\bf q}})$
In particular, this computation shows that if 
$f,g\in{\mathcal S}(U)$, we have  $(fg)_\H=f_\H g_\H= g_\H f_\H $,
so the map $f\mapsto f_\H$ is multiplicative. It is also clearly unital.

\section{The Quaternionic Cauchy Transform}

Using the $\M$-valued Cauchy kernel, we may define a concept of 
Cauchy transform, whose main properties will be discussed in this section.

The frequent use of versions of the Cauchy formula is simplified by adopting the following definition. Let $U\subset\C$ be open. An open subset 
$\Delta\subset U$ will be called a {\it Cauchy domain} (in $U$) if 
$\Delta\subset\bar{\Delta}\subset U$ and the boundary $\partial\Delta$ of $\Delta$ consists of a finite family of closed curves, piecewise smooth, positively oriented.
Note that a Cauchy domain is bounded but not necessarily connected.
\medskip

For a given open set $U\subset\C$, we denote by $\mathcal{O}(U,\M)$ the complex algebra  of all $\M$-valued analytic functions on $U$. 

If $U\subset\C$ is open and conjugate symmetric, let 
$\mathcal{O}_s(U,\M)$ be the real subalgebra of $\mathcal{O}(U,\M)$ consisting of all stem functions from $\mathcal{O}(U,\M)$.

Because $\C\subset\M$, we have $\mathcal{O}(U)\subset\mathcal{O}(U,\M)$, where $\mathcal{O}(U)$ is the complex algebra  of all 
complex-valued analytic functions on the open set $U$. Similarly, when $U\subset\C$ is open and conjugate symmetric, $\mathcal{O}_s(U)\subset\mathcal{O}_s(U,\M)$, where $\mathcal{O}_s(U)$ is the 
real subalgebra consisting of all functions $f$ from $\mathcal{O}(U)$  which are stem functions. 

As un example, if $\Delta\subset\C$ is an open disk centered at $0$, each function $F\in\mathcal{O}_s(\Delta,\M)$ can be represented 
as a convergent series $F(\zeta)=\sum_{k\ge 0}a_k\zeta^k,\,\zeta\in\Delta$, with $a_k\in\H$ for all $k\ge 0$. 


\begin{Def}\label{vect_fc}\rm Let $U\subset\C$ be a conjugate symmetric open set,  and let $F\in\mathcal{O}(U,\M)$. For every ${\bf q}\in U_\H$ we set 
\begin{equation}\label{Cauchy_vect}
C[F]({\bf q})=\frac{1}{2\pi i}\int_\Gamma F(\zeta)(\zeta-{\bf q})^{-1}d\zeta,
\end{equation}
where $\Gamma$ is the boundary of a Cauchy domain in $U$ containing the spectrum $\sigma({\bf q})$. The function $C[F]:U_\H\mapsto\M$
is called the {\it (quaternionic) Cauchy transform} of the function $F\in\mathcal{O}(U,\M)$. Clearly, the function $C[F]$ does not depend on the choice of $\Gamma$ because the function $U\setminus\sigma({\bf q})\ni\zeta\mapsto 
F(\zeta)(\zeta-{\bf q})^{-1}\in\M$ is analytic. 

We  put
\begin{equation}\label{imagqCt}
\mathcal{R}(U_\H,\M)=\{C[F];F\in\mathcal{O}(U,\M)\}.
\end{equation}
\end{Def} 


\begin{Pro}\label{right_reg} Let $U\subset\C$ be open and conjugate symmetric, and let $F\in\mathcal{O}(U,\M)$. Then  function  $C[F]\in \mathcal{R}(U_\H,\M)$ is slice regular on $U_\H$.
\end{Pro}

{\it Proof.}\, Let $F\in\mathcal{O}(U,\M)$, let
${\bf q}\in U_\H$ and let  $\Delta\ni\sigma({\bf q})$ be a conjugate symmetric Cauchy domain in $U$, whose boundary is denoted by $\Gamma$. We use the representation of $C[F]({\bf q})$
given by (\ref{Cauchy_vect}). Because we have
$$
\bar{\partial}_\mathfrak{s}((\zeta-{\bf q})^{-1})=\bar{\partial}_\mathfrak{s}((\zeta-x-y\mathfrak{s})^{-1})=0
$$
for ${\bf q}=x+y\mathfrak{s}\in\Delta_\H\cap(\R+\R\mathfrak{s})$,
via Example \ref{SliceregCk}, we infer that 
$$
\bar{\partial}_\mathfrak{s}(C[F]({\bf q}))=\frac{1}{2\pi i}\int_\Gamma F(\zeta)\bar{\partial}_\mathfrak{s}((\zeta-{\bf q})^{-1})
d\zeta=0,
$$
which implies the assertion.
\medskip

\begin{Rem}\rm (1) Because the function $F$ does not necessarily commute with the left multiplication by 
$\mathfrak{s}$, the choice of the right multiplication in the 
slice regularity is necessary to get the stated property of 
$C[F]$.

(2)  Let $r>0$ and let $U\supset\{\zeta\in\C;\vert\zeta\vert\le r\}$ be a conjugate symmetric open set. Then for every
$F\in\mathcal{O}(U,\M)$ one has
$$
C[F]({\bf q})=\sum_{n\ge0}\frac{F^{(n)}(0)}{n!} {\bf q}^n,\,\,
{\bf q}\in U_\H,\,\,
\Vert \bf q\Vert<r,
$$
where the series is absolutely convergent. Of course, using the
convergent series $(\zeta-{\bf q})^{-1}=
\sum_{n\ge0}\zeta^{-n-1}{\bf q}^n$ in $\{\zeta;\vert\zeta\vert=r\}$,
the assertion follows easily, via formula (\ref{Cauchy_vect}). Moreover, by Proposition \ref{right_reg}, the function $C[F]$ is a slice regular $\M$-valued function in $U_\H$. Nevertheless, we are particularly interested in slice regular $\H$-valued functions.
\end{Rem}


\begin{Thm}\label{vect_afc} Let $U\subset\C$ be a conjugate symmetric open set 
and let $F\in\mathcal{O}(U,\M)$. The Cauchy transform $C[F]$ is $\H$-valued if and only if $F\in\mathcal{O}_s(U,\M)$.
\end{Thm} 

{\it Proof.\,} We first fix a ${\bf q}\in U\setminus\R$.  If 
$\sigma({\bf q})=\{s_+,s_-\}$, the points $s_+,s_-$ are distinct and  not real. We then fix an $r>0$ sufficiently small such that,
setting $D_\pm:=\{\zeta\in U;\vert\zeta-s_\pm\vert\le r\}$, we
have 
$D_\pm\subset U$ and $D_+\cap D_-=\emptyset$. Then 
$$
C[F]({\bf q})=\frac{1}{2\pi i}\int_{\Gamma_+} F(\zeta)(\zeta-{\bf q})^{-1}d\zeta+\frac{1}{2\pi i}\int_{\Gamma_-} F(\zeta)(\zeta-{\bf q})^{-1}d\zeta,
$$
where $\Gamma_\pm$ is the boundary of $D_\pm$. We may write 
$F(\zeta)=\sum_{k\ge0}(\zeta-s_+)^ka_k$ with $\zeta\in D_+$,
$a_k\in\M$ for all $k\ge0$, as a uniformly convergent series.
Similarly, $F(\zeta)=\sum_{k\ge0}(\zeta-s_-)^kb_k$ with $\zeta\in D_-$, $b_k\in\M$ for all $k\ge0$, as a uniformly convergent series. 
 
Note that
$$
\frac{1}{2\pi i}\int_{\Gamma_+} F(\zeta)(\zeta-{\bf q})^{-1}d\zeta=\sum_{k\ge0}\left(a_k\frac{1}{2\pi i}\int_{\Gamma_+} (\zeta-s_+)^k(\zeta-{\bf q})^{-1}d\zeta\right)=
a_0\iota_+(\mathfrak{s}_{\tilde{\bf q}})
$$
because  we have
$$
\frac{1}{2\pi i}\int_{\Gamma_+} (\zeta-s_+)^k(\zeta-{\bf q})^{-1}d\zeta=({\bf q}-s_+)^k\iota_+(\mathfrak{s}_{\tilde{\bf q}})
$$
by the analytic functional calculus of ${\bf q}$ (see also
Lemma \ref{spproj}), which is equal
to $\iota_+(\mathfrak{s}_{\tilde{\bf q}})$ when $k=0$,  and it is equal to $0$ when $k\ge 1$, via the equality ${\bf q}\iota_+(\mathfrak{s}_{\tilde{\bf q}})=s_+\iota_+(\mathfrak{s}_{\tilde{\bf q}})$ 

\medskip

Similarly
$$
\frac{1}{2\pi i}\int_{\Gamma_-} F(\zeta)(\zeta-q)^{-1}d\zeta=\sum_{k\ge0}\left(b_k\frac{1}{2\pi i}\int_{\Gamma_-} (\zeta-s_-)^k(\zeta-{\bf q})^{-1}d\zeta\right)=b_0\iota_-(\mathfrak{s}_{\tilde{\bf q}})
$$
because, as above, we have
$$
\frac{1}{2\pi i}\int_{\Gamma_-} (\zeta-s_-)^k
(\zeta-{\bf q})^{-1})d\zeta=({\bf q}-s_-)^k\iota_-(\mathfrak{s}_{\tilde{\bf q}}),
$$
which is equal $\iota_-(\mathfrak{s}_{\tilde{\bf q}})$ when $k=0$, and it is equal to $0$ when $k\ge 1$. Consequently,
$$
C[F]({\bf q})=F(s_+)\iota_+(\mathfrak{s}_{\tilde{\bf q}})+F(s_-)\iota_-(\mathfrak{s}_{\tilde{\bf q}}),
$$
and the right hand side of this equality coincides with the expression from formula (5). 
 
Assume now that
$\sigma({\bf q})=\{s\}$,  where  $s:=s_+=s_-\in\R$. We fix an $r>0$ such that the set 
$D:=\{\zeta\in U;\vert\zeta-s\vert\le r\}\subset U$, 
whose boundary is denoted by $\Gamma$.  Then we have
$$
C[F]({\bf q})=\frac{1}{2\pi i}\int_{\Gamma} F(\zeta)(\zeta-{\bf q})^{-1}d\zeta=F(s),
$$
via the usual analytic functional calculus.

 In all of these situations, the element $C[F]({\bf q}))$ is equal to  the right hand side of formula (\ref{GFC2}). Therefore, we must have 
$C[F]({\bf q})\in\H$ if and only if $F(s_+)=\overline{F(s_-)}$, via Theorem \ref{sym_spec}. Consequently, $C[F]({\bf q})\in\H$ for all ${\bf q}\in U_\H$ if and only if $F:U\mapsto \M$ is a stem function.  
 

\begin{Rem}\label{bi_not}\rm (1) It follows from the proof of the previous theorem 
that the element $C[F]({\bf q}))$, given by formula (\ref{Cauchy_vect}), coincides with the element $F_\H({\bf q}))$ given by (\ref{GFC2}). To unify the notation, from now on this
element will be denoted by $F_\H({\bf q})$, whenever $F$ is a stem function, analytic or not.

(2) An important particular case is when let $f:U\mapsto\C$ is an analytic function, where  $U\subset\C$ is a conjugate symmetric open set. In this case we may also consider the (quaternionic) Cauchy transform of $f$ given by

\begin{equation}\label{Cauchy_vect1}
C[f]({\bf q})=\frac{1}{2\pi i}\int_\Gamma f(\zeta)(\zeta-{\bf q})^{-1}d\zeta,
\end{equation}
where  $\Gamma$ is the boundary of a Cauchy domain in $U$ containing the spectrum $\sigma({\bf q})$. According to Theorem 
\ref{vect_fc},  we have  $C[f]({\bf q})\in\H$ if and only if
$f(s_+({\bf q}))=\overline{f(s_-({\bf q}))}$ for all ${\bf q}\in U_\H$. In ohter words, $C[f]({\bf q})\in\H$ for all ${\bf q}\in\H$
if and only if $f$ is a stem function, that is 
$f\in\mathcal{O}_s(U)$. Of course, in this case 
 we may (and shall) also use the notation $C[f]=f_\H$. 
\end{Rem}

\section{Analytic Functional Calculus in Quaternionic Framework}
\label{AFCQ}

Let $\Omega\subset\H$ be a spectrally saturated open set, and 
let $U=\mathfrak{S}(\Omega)\subset\C$ (which is also open by
Remark \ref{consymr}(2)). We set
$$
\mathcal{R}(\Omega)=\{f_\H; f\in\mathcal{O}_s(U)\},
$$
which is an $\R$-algebra, and
$$
\mathcal{R}(\Omega,\H)=\{F_\H; F\in\mathcal{O}_s(U,\M)\},
$$
which, according to the next theorem, is a right $\mathcal{R}(\Omega)$-module. 

In fact, these $\R$-linear spaces have some important 
properties, as already noticed in a version of the next theorem (see Theorem 2 in \cite{Vas6}).


\begin{Thm}\label{H_afc} Let $\Omega\subset\H$ be a spectrally
saturated open set, and let $U=\mathfrak{S}(\Omega)$.
The space $\mathcal{R}(\Omega)$ is a unital commutative $\R$-algebra, the space $\mathcal{R}(\Omega,\H)$ is a right 
$\mathcal{R}(\Omega)$-module, the map
$$
{\mathcal O}_s(U,\M)\ni F\mapsto F_\H\in\mathcal{R}(\Omega,\H)
$$ 
is a right module isomorphism, and  its restriction
$$
{\mathcal O}_s(U)\ni f\mapsto f_\H\in\mathcal{R}(\Omega)
$$ 
is an $\R$-algebra isomorphism.

Moreover, for every polynomial\,\,$P(\zeta)=\sum_{n=0}^m a_n\zeta^n,\,\zeta\in\C$, with $a_n\in\H$  for all $n=0,1,\ldots,m$, we have  $P_\H(q)=\sum_{n=0}^m a_n q^n\in\H$ for all $q\in\H$.   
\end{Thm}

{\it Proof.\,} Thanks to Theorem \ref{vect_afc}, this statement 
is a particular case of Theorem \ref{gfc}. Indeed, the 
 $\R$-linear maps
$$
{\mathcal O}_s(U,\M)\ni F\mapsto F_\H\in\mathcal{R}(\Omega,\H),\,
{\mathcal O}_s(U)\ni f\mapsto f_\H\in\mathcal{R}(\Omega),\,\,
$$
are restrictions of the maps
$$
{\mathcal S}(U,\M)\ni F\mapsto F_\H\in\mathcal{F}(\Omega,\H),\,\,
{\mathcal S}(U)\ni f\mapsto f_\H\in\mathcal{F}(\Omega,\H),
$$ 
respectively. Moreover, they are $\R$-isomorphisms, the latter being actually unital and  multiplicative.   
Note that, in particular,  for every polynomial $P(\zeta)=\sum_{n=0}^m a_n\zeta^n$ with $a_n\in\H$  for all $n=0,1,\ldots,m$, we have  $P_\H({\bf q})=\sum_{n=0}^m a_n {\bf q}^n\in\H$ for all ${\bf q}\in\H$.


\begin{Rem}\label{n-deriv}\rm For every function $F\in\mathcal{O}_s(U,\M)$, the derivatives $F^{(n)}$ also belong to 
$\mathcal{O}_s(U,\M)$, where $U\subset\C$ is a conjugate symmetric open set.  

Now fixing $F\in\mathcal{O}_s(U,\M)$, we may define its 
{\it extended derivatives} with respect to the quaternionic variable via the formula

\begin{equation}\label{Cauchy_vect_deriv}
F^{(n)}_\H({\bf q})=\frac{1}{2\pi i}\int_\Gamma F^{(n)}(\zeta)(\zeta-{\bf q})^{-1}d\zeta,
\end{equation}
for the boundary $\Gamma$ of a Cauchy domain $\Delta\subset U$, $n\ge0$ an arbitrary integer, and $\sigma({\bf q})\subset\Delta$. 

In particular, if $\Delta$ is a disk centered at zero and 
$F\in\mathcal{O}_s(\Delta,\M)$, so we have a representation of
$F$   as a convergent series $\sum_{k\ge0}a_k\zeta^k$ with coefficients in $\H$, then 
(\ref{Cauchy_vect_deriv}) gives the equality
$F'_\H({\bf q})=\sum_{k\ge1}ka_k{\bf q}^{k-1}$, which looks like a (formal) derivative of the function $F_\H({\bf q})=\sum_{k\ge0}a_k{\bf q}^{k}$. 
\end{Rem}

The functions from the space $\mathcal{R}(\Omega,\H)$ admit a series
development around any real point of their domain of definition. 
In this sense, we have the following.


\begin{Pro}\label{series} Let $U\subset\C$ be conjugate symmetric, let $F\in\mathcal{O}_s(U,\M)$,  let $s_0\in\R\cap U$, and let $r>0$ be such that $D_0=\{\zeta\in U;\vert\zeta-s_0\vert< r\}\subset\bar{D_0}\subset U$. Then we have
$$
F_\H({\bf q})=\sum_{n\ge 0}\frac{F^{(n)}(s_0)}{n!}({\bf q}-s_0)^n,\,\,
\sigma({\bf q})\subset D_0.
$$
\end{Pro}

{\it Proof}\, Fixing ${\bf q}=x+y\mathfrak{s}$ with $x,y\in\R$
and $\mathfrak{s}\in\S$, such that $\sigma({\bf q})=\{x\pm iy\}\subset D_0$,
we must have
$$
\Vert{\bf q}-s_0\Vert^2=\vert x-s_0\vert^2+y^2=\vert x\pm iy-s_0\vert^2<r^2,
$$
implying the convergence of the series $\sum_{n\ge 0}r^{-n-1}({\bf q}-s_0)^n.$ Therefore,
$$
F_\H({\bf q})=\frac{1}{2\pi i}\int_{\Gamma_0}F(\zeta)(\zeta-{\bf q})^{-1}d\zeta=
 \frac{1}{2\pi i}\int_{\Gamma_0}F(\zeta)\sum_{n\ge 0}\frac{({\bf q}-s_0)^n}{(\zeta-s_0)^{n+1}}d\zeta=
$$
$$
\sum_{n\ge 0} \frac{F^{(n)}(s_0)}{n!}
({\bf q}-s_0)^n,
$$
where $\Gamma_0$ is the boundary of $D_0$.


\begin{Rem}\label{deriv}\rm As already noticed in the framework of 
\cite{Vas6}, Theorem \ref{H_afc} suggests 
a definition for $\H$-valued "analytic functions`` as elements of the set $\mathcal{R}(\Omega,\H)$, where $\Omega$ is a 
spectrally saturated open subset of $\H$. 
Because the expression "analytic function`` is quite improper 
in this context, the elements of $\mathcal{R}(\Omega,\H)$ will be called {\it Q-regular functions} on $\Omega$. In fact, 
the functions from $\mathcal{R}(\Omega,\H)$ are  Cauchy transforms of the stem functions from $\mathcal{O}_s(U,\M)$, with $U=\mathfrak{S}(\Omega)$. Moreover, as Proposition \ref{series} seems to suggest,
the functions from  $\mathcal{R}(\Omega,\H)$ are ''real analytic``
rather than ''analytic``.

Except for Theorem \ref{H_afc}, many other properties of $Q$-regular 
functions can be obtained directly from the definition, by 
recapturing the corresponding results from \cite{Vas6}. We omit the details.
\end{Rem}


\begin{Rem}\label{repres_form}\rm Let $U\subset\C$ be conjugate symmetric, let $x,y\in\R$ with $y\neq0$ and $z_\pm=x\pm iy\in U$, let  $F\in\mathcal{O}_s(U,\M)$, and let $\mathfrak{s}\in\S$

Assuming $y>0$, we consider the quaternions 
${\bf q}_\pm=x\pm y \mathfrak{s}$ for which 
for which
$s_+({\bf q}_\pm)=x+iy,\,s_-({\bf q}_\pm)=x-iy$. 

As we have
$\tilde{\bf q}_\pm=\pm y\mathfrak{s}$, then  
$s_{\tilde{\bf q}_\pm}=\pm\mathfrak{s}$, and $\iota_\pm(s_{\tilde{\bf q}_+})=(1\mp i\mathfrak{s})/2, 
\iota_\pm(s_{\tilde{\bf q}_-})=(1\pm i\mathfrak{s})/2.$ Therefore,
$$
2F_\H({\bf q}_+)=F(z_+)(1-i\mathfrak{s})+F(z_-)(1+i\mathfrak{s}),
$$
$$
2F_\H({\bf q}_-)=F(z_+)((1+i\mathfrak{s})+F(z_-)(1-i\mathfrak{s}).
$$
From these equations we deduce that
\begin{equation}\label{reprez1}
2F(z_+)=F_\H({\bf q}_+)(1-i\mathfrak{s})+F_\H({\bf q}_-)(1+i\mathfrak{s}),
\end{equation}
\begin{equation}\label{reprez2}
2F(z_-)=F_\H({\bf q}_+)(1+i\mathfrak{s})+F_\H({\bf q}_-)(1-i\mathfrak{s}).
\end{equation}

If $y<0$, for the quaternions ${\bf q}_\pm=x\pm y \mathfrak{s}$
we have $s_+({\bf q}_\pm)=x-iy,\,s_-({\bf q}_\pm)=x+iy$.
Moreover, as $\tilde{\bf q}_\pm=\pm y\mathfrak{s}$, then  
$s_{\tilde{\bf q}_\pm}=\mp\mathfrak{s}$, and $\iota_\pm(s_{\tilde{\bf q}_+})=(1\pm i\mathfrak{s})/2,\break 
\iota_\pm(s_{\tilde{\bf q}_-})=(1\mp i\mathfrak{s})/2.$
Therefore
$$
2F_\H({\bf q}_+)=F(z_-)(1+i\mathfrak{s})+F(z_+)(1-i\mathfrak{s}),
$$
$$
2F_\H({\bf q}_-)=F(z_-)((1-i\mathfrak{s})+F(z_+)(1+i\mathfrak{s}).
$$
These formulas lead again to equations (\ref{reprez1}) and 
(\ref{reprez2}). Consequently, we have the following.


\begin{Pro} Let $U\subset\C$ be conjugate symmetric, let $x,y\in\R$ with $x\pm iy\in U$, let $\mathfrak{s}\in\S$, and let $F\in\mathcal{O}_s(U,\M)$. Then we have the formulas
\begin{equation}\label{reprez}
F(x\pm iy)=F_\H(x\pm y\mathfrak{s})\left(\frac{1\mp i\mathfrak{s}}{2}\right)+
F_\H(x\mp y\mathfrak{s})\left(\frac{1\pm i\mathfrak{s}}{2}\right).
\end{equation}
\end{Pro} 

As the proof has been previously done, we only note that 
equality (\ref{reprez}) also holds  for $y=0$.
\end{Rem}


\begin{Lem}\label{equiv-ons-dom} Let $U\subset\H$ be a conjugate symmetric open set, let $\mathfrak{s}\in\S$ be fixed,
and let $\Psi:U_\mathfrak{s}\mapsto\H$ be  such that $\bar{\partial}_{\pm\mathfrak{s}}\Psi=0$. Then there exists  a function $\Phi\in\mathcal{R}(U_\H,\H)$ with $\Psi=\Phi\vert U_\mathfrak{s}$, where $U_\mathfrak{s}=\{x+y\mathfrak{s},x+iy\in U\}$.    
\end{Lem}

{\it Proof.}\, For arbitrary points $z_\pm=x\pm iy\in U$  with $x,y(\neq0)\in 
\R$, as in Remark \ref{repres_form}, we consider the quaternions ${\bf q}_\pm=x\pm y \mathfrak{s}$, so $s_+( {\bf q}_\pm)=x+i\vert y\vert$, and $s_( {\bf q}_\pm)=x-i\vert y\vert$. Inspired by formula (\ref{reprez}),  we set
$$
2F(z_+)=\Psi({\bf q}_+)(1-i\mathfrak{s})+\Psi({\bf q}_-)(1+i\mathfrak{s}),
$$
$$
2F(z_-)=\Psi({\bf q}_+)(1+i\mathfrak{s})+\Psi({\bf q}_-)(1-i\mathfrak{s}).
$$
 
Then we have
$$
2\frac{\partial F(z_+)}{\partial x}=\frac{\partial\Psi({\bf q}_+) }{\partial x}(1-i\mathfrak{s})+\frac{\partial\Psi({\bf q}_-) }
{\partial x}(1+i\mathfrak{s}),
$$
and
$$
2i\frac{\partial F(z_+)}{\partial y}=\frac{\partial\Psi({\bf q}_+) }{\partial y}\mathfrak{s}(1-i\mathfrak{s})+\frac{\partial\Psi({\bf q}_-)}{\partial y}(-\mathfrak{s})(1+i\mathfrak{s}),
$$
because $i(1-i\mathfrak{s})=\mathfrak{s}(1-i\mathfrak{s})$
and $i(1+i\mathfrak{s})=-\mathfrak{s}(1+i\mathfrak{s})$.

Therefore,
$$
\frac{\partial F(z_+)}{\partial x}+i\frac{\partial F(z_+)}{\partial y}=(\bar{\partial}_\mathfrak{s}\Psi({\bf q}_+))(1-i\mathfrak{s})+(\bar{\partial}_\mathfrak{-s}\Psi({\bf q}_-))(1+i\mathfrak{s})=0,
$$
showing that the function $z_+\mapsto F(z_+)$ is analytic in $U$. 

Because $\overline{F(z_-)}=\overline{F(\overline{z_+})}=F(z_+)$, and when $y=0$ we have
$\overline{F(z_-)}=F(z_+)=F(x)$, we have constructed a function
$F\in\mathcal{O}_s(U,\M)$. Hence, taking $\Phi=F_\H$, we have
$\Phi\in\mathcal{R}(U_\H,\H)$ with $\Psi=\Phi\vert U_\mathfrak{s}$, via Remark \ref{repres_form}.

 
\begin{Thm}\label{equiv-ons-dom1} Let $\Omega\subset\H$ be a spectrally saturated open set, and let $\Phi:\Omega\mapsto\H$.
 The following conditions are equivalent:

$(i)$ $\Phi$ is  a slice regular function;  

$(ii)$ $\Phi\in\mathcal{R}(\Omega,\H)$, that is, $\Phi$ is 
$Q$-regular. 
\end{Thm}

{\it Proof.}\, If $\Phi\in\mathcal{R}(\Omega,\H)$, then $\Phi$
is slice regular, by Lemma \ref{right_reg}, so $(ii)\Rightarrow(i)$.

Conversely, let $\Phi$ be slice regular in $\Omega$. Fixing a 
$\mathfrak{s}\in\S$, we have $\bar{\partial}_{\pm\mathfrak{s}}\Phi_{\mathfrak{s}}=0$, where $\Phi_{\mathfrak{s}}=\Phi\vert U_{\mathfrak{s}}$. It follows from Lemma \ref{equiv-ons-dom}
that there exists $\Psi\in\mathcal{R}(U_\H,\H)$ with
$\Psi_{\mathfrak{s}}=\Phi_{\mathfrak{s}}$. This implies that 
$\Phi=\Psi$, because both $\Phi,\Psi$ are uniquely determined by 
$\Phi_{\mathfrak{s}}, \Psi_{\mathfrak{s}}$, respectively, the former by (the right hand version of) Lemma 4.3.8 in \cite{CoSaSt}, and the latter by  Remark \ref{zeros}. Consequently, we also have $(i)\Rightarrow(ii)$.
\medskip

\noindent{\bf Final Remark} Theorem \ref{equiv-ons-dom1} allows us 
to obtain the properties of what we called $Q$-regular functions via 
those of the slice regular functions, as in \cite{CoSaSt}. Clearly, there is no coincidence that a  result like Proposition \ref{series} looks like Corollary 4.2.3 from \cite{CoSaSt}.  Nevertheless, as mentioned before, they can also be obtained directly, with the our techniques (see also \cite{Vas6}).

\end{document}